\documentclass[11pt]{amsart}
\usepackage{amsfonts}
\usepackage{graphicx}
\usepackage{amsmath}
\usepackage{amssymb}
\newtheorem{thm}{Theorem}[section]

\newtheorem{prop}[thm]{Proposition}
\theoremstyle{definition}

\newtheorem{rem}[thm]{Remark}
\numberwithin{equation}{section}

\begin{document}
\title[]{A perturbed time-dependent second-order sweeping process}%
\author{Fatine. ALIOUANE}
\address{F. Aliouane, Laboratoire LAOTI, Universit\'e Mohammed Seddik Benyahia de Jijel}
 \email{faliouane@gmail.com}
 \author{Dalila. AZZAM-LAOUIR}
 \address{D. Azzam-Laouir, Laboratoire LAOTI, Universit\'e Mohammed Seddik Benyahia de Jijel}
 \email{laouir.dalila@gmail.com, dalilalaouir@univ-jijel.dz}
\begin{abstract}In this paper we prove the existence of solutions for a second order sweeping
process with a Lipschitz single valued perturbation by transforming
it to a first order problem.
\end{abstract}
\maketitle \vskip4mm \noindent{\footnotesize\textit{Keywords:}
Lipschitz perturbation; normal cone; sweeping process}\vskip4mm
\noindent{\footnotesize\textit{2000 Mathematics Subject
Classification:} 34A60, 28A25, 28C20}\vskip6mm
\section{\bf Introduction}
In \cite{A}, we proved the existence of solutions of the following
second order sweeping process
$$(\mathcal{P}_f)
\begin{cases}
   -\ddot{x}(t)\in
   N_{K(t)}(\dot{x}(t))+f(t,x(t),\dot x(t)),\;a.e.\;t\in I;\\
   \dot{x}(t)\in K(t),\;\forall t\in I;\\
   x(0)=x_0;\;\dot{x}(0)=u_0.
\end{cases}$$
In the proof, we used a catching-up algorithm, that is, we
discretized the time interval $I=[0,T]$ on subintervals
$I_{n,i}=(t_{n,i},t_{n,i+1}]$ such that $t_{n,i}=ih$ with
$h=\frac{T}{n}$ for all $0\leq i\leq n$.

On such interval $I_{n,i}$, we  defined two discrete sequences
$(u_{n,i})$, $(x_{n,i})$ by
$$u_{n,i+1}=P_{K(t_{n,i+1})}\bigg(u_{n,i}-\int^{t_{n,i+1}}_{t_{n,i}}f(s,x_n({t_{n,i}}),u_{n,i})ds\bigg)$$
and $$ x_{n,i+1}=x_{n,i}+hu_{n,i}.$$ This algorithm is well defined
because the sets $K(t_{n,i})$ were supposed to be r-prox-regular.
Through these discrete sequences, we constructed two approximating
ones, $(u_n)$ and $(x_n)$, and we proved that they converge to the
desired solution. To this end, a compactness assumption on the sets
$K(t)$ was inevitable. This method was introduced in the 70's by
Moreau in \cite{M}, to solve a first order sweeping process, in the
setting where all the sets were assumed to be convex. Since then, a
lot of work has been investigated in this area and many authors
studied the problem of sweeping process   in different cases, for
the first order we can see for example \cite{adly-hadd-thib, A.I.T,
{A.M.T}, {B.T}, {C.D.V}, {C.I.Y}, {C.M}, {C.S.T}, {C.G},
{Chem-Mont}, {DE}, {H.T}, {Kunz-Mont}, {t}, {Vala}} and the
references therein and for the second order sweeping processes we
can refer to \cite{A.K, {A.N}, {A1}, {Az}, {A.I}, {B}, {B.A}, {C},
{Mar}} and the references therein.

Our aim in this paper, is to  give an  existence result for the
 second order sweeping process $(\mathcal{P}_f)$,
 without any compactness assumption on $K(t)$. The idea is to transform this second order problem to a
 first order sweeping process and use the following result established  by Edmond and
Thibault in \cite{ED}.
\begin{thm}\label{t1}
Let $H$ be a Hilbert space and $C:I=[T_0,T]\rightrightarrows H$ be
such
that\\
(\textbf{$H_1$}) for each $t\in I$, $C(t)$ is  a nonempty closed subset of $H$ which is $r$-prox-regular.\\
(\textbf{$H_2$}) $C(t)$ varies in an absolutely continuous way, that
is, there exists an absolutely continuous function
$v(.):I\rightarrow \mathbb{R}$ such that, for any $y\in H$ and
$s,t\in I,$
$$|d(y,C(t))-d(y,C(s))|\leq |v(t)-v(s)|.$$
Let $f:I\times H
\longrightarrow H$ be a separately measurable map on $I$ such that \\
\textbf{(i)} for every $\eta>0,$ there is a non-negative function
$\gamma_\eta(.)\in \mathbf{L}^1(I,\mathbb{R})$ such that for all
$t\in I$ and for any $(x,y)\in
\overline{\mathbf{B}}(0,\eta)\times\overline{\mathbf{B}}(0,\eta),$
$$\|f(t,x)-f(t,y)\|\leq \gamma_\eta(t)\|x-y\|;$$
\textbf{(ii)} there exists a non-negative function $\beta(.)\in
\mathbf{L}^1(I,\mathbb{R})$ such that, for all $t\in I$ and for all
$x\in \bigcup_{s\in I}C(s)$,
$$\|f(t,x)\|\leq\beta(t)(1+\|x\|).$$ Then, for any $x_0\in
C(T_0)$ the following perturbed sweeping process
$$(\textbf{SPP})
\begin{cases}
   -\dot{u}(t)\in N_{C(t)}(u(t))+f(t,x(t)),\;a.e.\;t\in I,\\
   u(T_0)=x_0
\end{cases}$$
has one and only one absolutely continuous solution $u(.)$. This
solution satisfies
$$\|\dot{u}(t)+f(t,u(t))\|\leq (1+l)\beta(t)+|\dot{v}(t)|\;\textmd{p.p.}\;t\in I$$
and$$\|f(t,u(t))\|\leq(1+l)\beta(t)\;\textmd{p.p.}\;t\in I,$$where
$$l:=\|x_0\|+exp\{2\int_{T_0}^T\beta(s)ds\}\int_{T_0}^T \big(2\beta(s)(1+\|x_0\|)+|\dot{v}(t)|\big)\,ds.$$
\end{thm}
\section{\bf Preliminaries}
In the sequel $H$ denotes a real  Hilbert space, $I=[0,T],\;T>0$. We
denote by $\mathbb{B}$ the closed unit ball and for $\eta>0$,
$\overline{\mathbf{B}}(0,\eta)$ is the closed ball of radius $\eta$
and center $0$.

We will denote by $\mathbf{C}(I,H)$ or $\mathbf{C}_H(I)$ the Banach
space of all continuous maps from $I$ to $H$ equipped with the norm
of the uniform convergence $\|.\|_{\mathbf{C}}$, and by
$\mathbf{C}^1_H(I)$ the Banach space of all continuous maps from $I$
to $H$ having continuous derivatives, equipped with the norm
$\|u\|_{\mathbf{C}^1}=\max\big(\|u\|_{\mathbf{C}}, \|\dot
u\|_{\mathbf{C}}\big)$. We denote by $\mathcal{L}(I)$ the
$\sigma$-algebra of Lebesgue measurable subsets of $I$.
$(\mathbf{L}_{H}^1(I),\|.\|_1)$ is the Banach space of
Lebesgue-Bochner integrable $H$-valued maps

We say that a map $u:[0,T]\rightarrow H$ is absolutely continuous if
there is a map $v\in\mathbf{L}_{H}^1(I) $ such that
$u(t)=u(0)+\int_{0}^tv(s)ds,$ for all $t\in I,$ in this case
$v=\dot{u}$ a.e.

 For a subset $A$  of $H$ and $x\in H$ we denote $d(x,A)=\inf_{y\in
 A}\|x-y\|$.

For a fixed $r>0$, the set $S$ is said to be $r$-prox-regular (or
uniformly prox-regular (see \cite{15}) or equivalently
$r$-proximally smooth (see \cite{24}) if and only if for all
$\bar{x}\in S$ and all $0\neq \xi\in N^P(S,x)$
$$S\cap int(\bar{x}+r\frac{\xi}{\|\xi\|}+r\mathbb{B})=\emptyset,\;\textmd{i.e.}\;,\langle\frac{\xi}{\|\xi\|},x-\bar{x}\rangle\leq \frac{1}{2r}\|x-\bar x\|^2,$$
for all $x\in S$. We make the convention $\frac{1}{r}=0$ for
$r=+\infty$. Recall that for $r=+\infty$, the $r$-prox-regularity of
$S$ is equivalent to its convexity.

Next, we give a useful Proposition needed in the proof of our
theorem which was communicated to us by Prof. Lionel. Thibault.
\begin{prop}\label{p1}
Let $S$, $S'$ two non empty closed sets of $H$ such that $S$ is
$r$-prox-regular and $S'$ is $r'$-prox-regular. Then, $S\times S'$
is $\rho$-prox-regular with $\rho=\min \{r,r'\}$.
\end{prop}

\proof For every $x,x'\in H$, we have $$d^2((x,x'),S\times
S')=d^2(x,S)+d^2(x',S').$$ Since $S$ is $r$-prox-regular then,
$d^2(.,S)$ is Fr\'echet-differentiable on \\$U_r(S)=\{x\in
H:\;d(x,S)<r\}$, and since $S'$ is $r'$-prox-regular then,
$d^2(.,S')$ is Fr\'echet-differentiable on $U_{r'}(S')=\{x\in
H:\;d(x',S')<r'\}$. We conclude that $d^2((.,.),S\times S')$ is
Fr\'echet-differentiable on $U_r(S)\times
U_{r'}(S')$.\\
Assume that $\rho=\min\{r,r'\}$. Put $$U_{\rho}(S\times
S')=\{(x,x')\in H\times H,\;d((x,x'),S\times S')<\rho\}$$ and let us
show that $$U_{\rho}(S\times S')\subset U_r(S)\times U_{r'}(S').$$
Let $(x,x')\in U_{\rho}(S\times S')$, so that, $d((x,x'),S\times
S')<\rho$, and by the last inequality, we have
$$\sqrt{d^2(x,S)+d^2(x',S')}<\rho=\min\{r,r'\}.$$ Therefore,
$$d(x,S)\leq \sqrt{d^2(x,S)+d^2(x',S')}<\min\{r,r'\}.$$ That is
 $d(x,S)<r$. By the same way
$d(x',S')<r'$. Then, $d^2((.,.),S\times S')$ is
Fr\'echet-differentiable on $U_{\rho}(S\times S')$. \hfill
$\blacksquare$

\section{\bf Main result}
\begin{thm}\label{th2}
Let $I=[0,T]\;\;(T>0)$ and $H$ be a real Hibert space. Let
$f:I\times H\times H \rightarrow H$ be a separately measurable map
on $I$ satisfying: \\ \textbf{(i)} for all $\eta>0,$ there exists a
non-negative function $k_\eta(.)\in \mathbf{L}^1_{\mathbb{R}}(I)$
such that for all $t\in I$ and for any $(x,y),(u,v)\in
\overline{\mathbf{B}}(0,\eta)\times\overline{\mathbf{B}}(0,\eta),$
$$\|f(t,x,u)-f(t,y,v)\|\leq k_\eta(t)(\|x-y\|+\|u-v\|);$$
\textbf{(ii)} for some non-negative real function  $c\in
\mathbf{L}^1_{\mathbb{R}}(I)$,
$$\|f(t,x,u)\|\leq c(t)(1+\|x\|+\|u\|),\;\textmd{for all}\;(t,x,u)\in
I\times H\times H.$$

Let $r>0$ and $K:I\rightrightarrows H$ be a set-valued map taking
nonempty, closed and r-prox-regular values. We assume that $K(.)$
moves in an absolutely continuous way, that is, there exists an
absolutely continuous real function $a(.)$ such that for all $t,s\in
I$ and $u\in H$
\begin{equation}\label{eq1}
|d(u, K(t))-d(u, K(s))|\leq |a(t)-a(s)|.
\end{equation}
 Let $x_0\in H$
and $u_0\in K(0)$.  Then the differential inclusion
$(\mathcal{P}_f)$ has an absolutely continuous solution $x(.)\in
\mathbf{C}_{H}^1(I)$.
\end{thm}
\proof

 Let us
consider the set-valued map $C:I\rightrightarrows H\times H$ defined
by  $C(t)=K(t)\times H$ for all $t\in I$, and  the mapping
$g:I\times H\times H\rightarrow H\times H$ defined by
$$g(t,x,y)=\big(f(t,y,x),-x\big),\;\forall (t,x,y)\in I\times H\times H.$$

We will show in the following, that the set-valued map $C$ is absolutely continuous.\\
Using the fact that $K(.)$ satisfies \eqref{eq1}, we have for all
$t,s\in I$ and $x,y\in H$,
\begin{align*}
\mid d\big((x,y),C(t)\big)-d\big((x,y),C(s)\big)\mid=&\mid
d\big((x,y),K(t)\times H\big)-d\big((x,y),K(s)\times
H\big)\mid\\
=&\mid d(x,K(t))+d(y,H)-d(x,K(s))-d(y,H)\mid\\=&\mid
d(x,K(t))-d(x,K(s))\mid\leq |a(t)-a(s)|.
\end{align*}
On the other hand, by virtue  of  Proposition \ref{p1}, it is clear
that $C$ has $\rho$-prox-regular values with $\rho=\min
\{r,+\infty\}=r.$ Then, assumptions ($H_1$) and $(H_2)$ of Theorem
\ref{t1} are fulfilled  on $C$.\\ Next, we will prove that $g$
satisfies the hypothesis of the same theorem. \\For all $\eta>0$,
$t\in I$ and all $(U,V)\in \overline{\mathbf{B}}_{H\times
H}(0,\eta)$ ($U=(u,u'),\;V=(v,v')$) we have,
\begin{align*}
\|g(t,U)-g(t,V)\|=&\|(f(t,u',u),-u)-(f(t,v',v),-v)\|\\
=&\|(f(t,u',u)-f(t,v',v),-u+v)\|\\
=&\|f(t,u',u)-f(t,v',v)\|+\|u-v\|\\
\leq&k_\eta(t)\big(\|u-v\|+\|u'-v'\|+\|u-v\|\big)\\
\leq&(k_\eta(t)+1)\big(\|u-v\|+\|u'-v'\|\big)\\
=&(k_\eta(t)+1)\|(u,u'),(v,v')\|=:\gamma_\eta(t)\|U-V\|.
\end{align*}
So that $\textbf{(i)}$ is clearly satisfied.  Concerning
$\textbf{(ii)}$, we have for every $t\in I$ and every $U\in
\displaystyle\bigcup_{s\in I}C(s)$ ($U=(u,u')$),
\begin{align*}
\|g(t,U)\|=&\|(f(t,u',u),-u)\|\\
=&\|f(t,u',u)\|+\|u\|\\
\leq& c(t)(1+\|u\|+\|u'\|)+\|u\|\\
\leq& (c(t)+1)(1+\|u\|+\|u'\|)\\
=& (c(t)+1)(1+\|U\|)=:\beta(t)(1+\|U\|).
\end{align*}
We conclude that all the hypothesis of Theorem \ref{t1} are
satisfied. Consequently, there exists a unique absolutely continuous
solution $X(.)=(y(.), x(.))$ of the problem $(\mathbf{SPP})$, that
is
$$
\begin{cases}
   -\dot{X}(t)\in
   N_{C(t)}(X(t))+g(t,X(t)),\;a.e.\;t\in I;\\
   X(t)\in C(t),\;\forall t\in I;\\
   X(0)=(u_0, x_0).
\end{cases}$$
This entails that for almost every $t\in I,$

$$-(\dot{y}(t),\dot{x}(t))\in N_{K(t)\times H}(y(t), x(t))+\big(f(t,x(t),y(t)),-y(t)\big),$$
then,
$$-(\dot{y}(t),\dot{x}(t))\in N_{K(t)}(y(t))\times \{0\}+\big(f(t,x(t),y(t)),-y(t)\big),$$
i.e.,
$$-\dot{y}(t)\in N_{K(t)}(y(t))+f(t,x(t),y(t)),$$
and $\dot x(t)=y(t)$. These relations mean that
$$-\ddot{x}(t)\in N_{K(t)}(\dot x(t))+f(t,x(t),\dot x(t)),\;\;a.e.\,t\in I.$$
Moreover, since for all $t\in I$, $X(t)\in C(t)$, we get
$\dot{x}(t)\in K(t)$ for all $t\in I$, and
$X(0)=(\dot{x}(0),x(0))=(u_0,x_0)$. Therefore $x(.)$ is a solution
in $\mathbf{C}_E^1(I)$ of our problem $(\mathcal{P}_f)$.
\hfill$\blacksquare$

\begin{rem}
This work was established in the PhD thesis of the first author
which was supervised by the second author of this paper. See
\cite{Fatine}.
\end{rem}

\end{document}